\documentclass[letterpaper, 10pt, conference]{ieeeconf}      

\IEEEoverridecommandlockouts                              
\newcommand{\ud}{\,\mathrm{d}}

\usepackage{graphics} 
\usepackage{epsfig} 
\usepackage{mathptmx} 
\usepackage{times} 
\usepackage{amsmath} 
\usepackage{amssymb}  
\usepackage{mathtools}
\usepackage{tabularx}
\usepackage{multirow}

\newtheorem{theorem}{Theorem}

\newtheorem{remark}{Remark}

\newtheorem{corollary}{Corollary}

\usepackage{algorithm}
\usepackage{algorithmic}

\newlength{\noteWidth}
\long\def\notes#1{\ifinner
             {\tiny #1}
             \else
              \marginpar{\parbox[t]{\noteWidth}{\raggedright\tiny #1}}
               \fi}

\def\notes#1{\typeout{#1 !!!}}  

\newcounter{rmnum}

\newenvironment{romannum}{\begin{list}{{\upshape (\roman{rmnum})}}{\usecounter{rmnum}
\setlength{\leftmargin}{6pt}
\setlength{\rightmargin}{4pt}
\setlength{\itemindent}{-1pt}
}}{\end{list}}

\newcounter{anum}

\def\Sec#1{Sec.~\ref{#1}}
\def\Fig#1{Fig.~\ref{#1}}
\def\Thm#1{Thm.~\ref{#1}}

\def\head#1{\noindent \textit{#1}\ }

\def\IEEEQEDclosed{\mbox{\rule[0pt]{1.3ex}{1.3ex}}}

\def\qed{\nobreak\hfill\IEEEQEDclosed}

\def\E{{\sf E}}
\def\P{{\sf P}}

\def\Re{\mathbb{R}}
\def\Se{\mathbb{S}}

\def\Inov{I}
\def\v{{\sf K}}
\def\Prob{{\sf P}}

\def\clL{{\cal L}}

\def\phatm{\rho_m}

\def\mth{m^{\text{th}}}

\def\UZ{\mathcal{Z}_t}
\def\Ldag{\mathcal{L}^\dagger}
\newcommand{\ith}{i^{\text{th}}}

\title{\LARGE \bf
Interacting Multiple Model-Feedback Particle Filter \\ for
Stochastic Hybrid Systems}

\author{Tao Yang,  Henk A. P. Blom, Prashant G. Mehta
\thanks{Financial support from the AFOSR grant FA9550-09-1-0190 and the NSF grant EECS-0925534 is gratefully acknowledged.}
\thanks{T. Yang and P. G. Mehta are with the Coordinated Science Laboratory and
the Department of Mechanical Science and Engineering at the University of
Illinois at Urbana-Champaign (UIUC)
{\tt\scriptsize taoyang1@illinois.edu; mehtapg@illinois.edu}}
\thanks{H. A. P. Blom is with National Aerospace Laboratory NLR and with Delft University of Technology, both in the Netherlands
{\tt\scriptsize blom@nlr.nl}}
}

\begin{document}

\maketitle
\thispagestyle{empty}

\begin{abstract}

In this paper, a novel feedback control-based particle
filter algorithm for the continuous-time stochastic hybrid system estimation
problem is presented. This particle filter is referred to as the
interacting multiple model-feedback particle filter (IMM-FPF), and is
based on the recently developed feedback particle filter~\cite{taoyang_cdc12,taoyang_cdc11,taoyang_acc11}.
The IMM-FPF is comprised of a series of parallel FPFs, one for
each discrete mode, and an exact filter recursion for the mode
association probability. The proposed IMM-FPF represents a
generalization of the Kalman-filter based IMM algorithm to
the general nonlinear filtering problem.

The remarkable conclusion of this paper is that the IMM-FPF algorithm retains
the innovation error-based feedback structure even for the
nonlinear problem. The interaction/merging process is also
handled via a control-based approach. The theoretical results
are illustrated with the aid of a numerical example problem
for a maneuvering target tracking application.
\end{abstract}

\section{INTRODUCTION}
\label{sec:intro}
State estimation for stochastic hybrid systems (SHS) is important to a
number of applications, including air and missile defense systems,
air traffic control, satellite surveillance, statistical pattern
recognition, remote sensing, autonomous navigation and robotics~\cite{BarshalomLiKirubarajan01}. A typical problem formulation involves estimation of
a partially observed stochastic process with both continuous-valued and discrete-valued states.

An example of the SHS estimation is the problem of tracking a maneuvering target
(hidden signal) with noisy radar measurements. In this case, the continuous-valued states are target positions and velocities, while the discrete-valued states represent the distinct dynamic model types (e.g., constant velocity or white noise acceleration model) of the target. The discrete signal model types are
referred to as {\em modes}. Since the time of target maneuvers is random,
there is model association uncertainty in the sense that one can not
assume, in an apriori fashion, a fixed dynamic model of the target.

Motivated in part by target tracking applications, we consider models of SHS where the continuous-valued state process is modeled using a stochastic differential equation (SDE), and the discrete-valued state process is modeled as a Markov chain. The estimation objective is to estimate (filter) the hidden states given noisy observations.

Given the number of applications, algorithms for SHS filtering problems
have been extensively studied in the past;
cf.,~\cite{BarshalomLiKirubarajan01,Mazor98} and references therein. A typical SHS filtering algorithm is comprised of three parts:
\begin{romannum}
\item A {\em filtering} algorithm to estimate the continuous-valued state given the mode,
\item An {\em association} algorithm to associate
  modes to signal dynamics,
\item A {\em merging process} to combine the results of i) and ii).
\end{romannum}

Prior to mid-1990s, the primary tool for filtering was a Kalman
filter or one of its extensions, e.g., extended Kalman filter.
The limitations of these tools in applications arise on account
of nonlinearities, not only in dynamic motion of targets (e.g.,
drag forces in ballistic targets) but also in the measurement
models (e.g., range or bearing).  The nonlinearities can lead
to a non-Gaussian multimodal conditional distribution.  For
such cases, Kalman and extended Kalman filters are known to
perform poorly; cf.,~\cite{Ristic_book_2004}. Since the advent
and wide-spread use of particle
filters~\cite{gorsalsmi93,DouFreGor01}, such filters are
becoming increasing relevant to SHS estimation for target
tracking applications; cf.,~\cite{Ristic_book_2004} and
references therein.

The other part is the mode association algorithm. The purpose of the
mode association algorithm is to determine the conditional
probability for the discrete modes.  

In a discrete-time setting, the exact solution to problems (i)-(iii) is given by a Multiple Model (MM) filter which has an (exponentially with time) increasing number of filters, one for each possible mode history. Practically, however, the number of filters has to be limited, which leads to the classical Generalised Pseudo-Bayes estimators of the first and second order (GPB1 and GPB2) and the Interacting Multiple Model (IMM) filter~\cite{BarshalomLiKirubarajan01}. For some SHS examples, however, it was already shown in~\cite{BlomBarshalom88TAC} that these low-dimensional filters do not always perform well. This has led to the development of two types of particle filters for SHS:
\begin{romannum}
\item The first approach is to apply the standard particle filtering approach to the joint continuous-valued state and discrete-valued mode process~\cite{Irwin00},\cite{Musso01}.
\item The second approach is to exploit Rao-Blackwellization, in the sense of applying particle filtering for the continuous-valued state, and exact filter recursions for the discrete-valued modes~\cite{Boers04},\cite{BlomBloem04},\cite{BlomBloem07}.
\end{romannum}

\smallskip

In this paper, we consider a continuous-time filtering problem for SHS and develop a novel feedback control-based particle filter algorithm, where the particles represent continuous-valued state components (case (ii)). The proposed algorithm is based on the feedback particle filter\;(FPF) concept introduced by
us in earlier
papers~\cite{taoyang_acc11},\cite{taoyang_cdc11},\cite{taoyang_cdc12}. A feedback
particle filter is a controlled system to approximate the
solution of the nonlinear filtering task. The filter has a
feedback structure similar to the Kalman filter: At each time
$t$, the control is obtained by using a proportional gain
feedback with respect to a certain modified form of the
innovation error.  The filter design amounts to design of the
proportional gain -- the solution is given by the Kalman gain
in the linear Gaussian case.

In the present paper, we extend the feedback particle filter to SHS estimation problems.  We refer to the resulting algorithm as the Interacting Multiple Model-Feedback Particle Filter (IMM-FPF).  As the name suggests, the proposed algorithm represents a generalization
of the Kalman filter-based IMM algorithm now to the general nonlinear filtering problem.

One remarkable conclusion of our paper is that the IMM-FPF
retains the innovation error-based feedback structure even for
the nonlinear problem. The interaction/merging process is also
handled via a control-based approach.  The innovation error-based
feedback structure is expected to be useful because of the
coupled nature of the filtering and the mode association
problem.  The theoretical results are illustrated with a numerical
example.

The outline of the remainder of this paper is as follows: The exact
filtering equations appear in~\Sec{sec:exact_eqn}.  The IMM-FPF is introduced in~\Sec{sec:Ifilter} and the numerical example is described in~\Sec{sec:numerics}.

\section{Problem formulation and exact filtering equations}
\label{sec:exact_eqn}
In this section, we formulate the continuous-time SHS filtering problem, introduce the notation, and summarize the exact filtering equations (see~\cite{Blom84,Loparo86TAC,Blom_cdc12} for standard references).
For pedagogical reason, we limit the considerations to scalar-valued signal and observation processes. The generalization to multivariable case is straightforward.

\subsection{Problem statement, Assumptions and Notation}
The following notation is adopted:
\begin{romannum}
\item At time $t$, the signal state is denoted by $X_t\in\Re$.
\item At time $t$, the mode random variable is denoted as $\theta_t$, defined on a state-space comprising of the standard basis in $\Re^M$: $\{e_1,e_2,\hdots,e_M\} =: \Se$. It associates a specific mode to the signal: $\theta_t = e_m$ signifies that the dynamics at time $t$ is described by the $m^{\text{th}}$ model.
\item At time $t$, there is only one observation $Z_t\in \Re$. The observation history~(filtration) is denoted as $\UZ:= \sigma(Z_s:s \leq t)$.
\end{romannum}

The following models are assumed for the three stochastic processes:
\begin{romannum}

\item The evolution of the continuous-valued state $X_t$ is described by a stochastic differential equation with discrete-valued coefficients:
\begin{equation}
\ud X_t = a(X_t,\theta_t) \ud t + \sigma(\theta_t) \ud B_t, \quad \label{eqn:process_model}
\end{equation}
where $B_t$ is a standard Wiener process. We denote $a^m(x) := a(x,e_m) $ and $\sigma^m := \sigma(e_m)$.
\item The discrete-valued state (mode) $\theta_t$ evolves as a Markov chain in continuous-time:
\begin{equation}
\P(\theta_{t+\delta}=e_l|\theta_t=e_m) = q_{ml}\delta + o(\delta),\quad m \neq l.\label{eqn:Markov_Process}
\end{equation}
The generator for this jump process is denoted by a stochastic matrix $Q$ whose $ml^{\text{th}}$ entry is $q_{ml}$ for $m\neq l$. The initial distribution is assumed to be given.

\item At time $t$, the observation model is given by,
\begin{equation}
\ud Z_t = h(X_t, \theta_t)\ud t + \ud W_t,\label{eqn:observ_model}
\end{equation}
where $W_t$ is a standard Wiener process assumed to be independent of $\{B_t\}$. We denote $h^m(x) := h(x,e_m)$.
\end{romannum}

The filtering problem is to obtain the posterior distribution of $X_t$ given $\UZ$. 

\subsection{Exact Filtering Equations}

The following distributions are of interest:
\begin{romannum}
\item $q_m^\ast(x,t)$ defines the joint conditional distribution of $(X_t,\theta_t)^T$ given $\mathcal{Z}_t$, i.e.,
\begin{equation*}
\int_{x\in A} q^\ast_m(x,t) \ud x = \P\{ [X_t \in A, \theta_t = e_m] | \UZ\}, 
\end{equation*}
for $A \in \mathcal{B}(\Re)$ and $m\in \{1,\hdots,M\}$. We denote $q^\ast(x,t) := (q_1^\ast(x,t),q_2^\ast(x,t),\hdots,q_M^\ast(x,t))^T$, interpreted as a column vector.
\item $p^\ast(x,t)$ defines the conditional dist. of $X_t$ given $\UZ$:
\begin{equation*}
\int_{x\in A} p^\ast(x,t) \ud x = \P\{X_t \in A | \UZ\},\quad A \in \mathcal{B}(\Re).
\end{equation*}
By definition, we have $p^\ast(x,t) = \sum_{m=1}^M q_m^\ast(x,t)$.
\item $\mu_t:= (\mu_t^1,\hdots,\mu_t^M)^T$ defines the probability mass function of $\theta_t$ given $\UZ$ where:
\begin{equation}
\mu_t^m = \P\{\theta_t=e_m|\UZ\},\quad m=1,\hdots, M.\label{eqn:def_mu_m}
\end{equation}
By definition $\mu_t^m = \int_{\Re} q_m^\ast(x,t)\ud x$. 
\item $\rho_m^\ast(x,t)$ defines the conditional dist. of $X_t$ given $\theta_t = e_m$ and $\UZ$. For $\mu_t^m \neq 0$:
\begin{equation}
\rho_m^\ast(x,t) := \frac{q_m^\ast(x,t)}{\mu_t^m}, \quad m = 1,\hdots,M, \label{eqn:def_rho_m}
\end{equation}
Denote $\rho^\ast (x,t) = (\rho_1^\ast(x,t),\hdots,\rho_M^\ast(x,t))^T$.

\end{romannum}

We introduce two more notations before presenting the exact filtering equations for these density functions:
\begin{romannum}
\item $\hat{h}_t := \E[h(X_t,\theta_t)|\UZ] = \sum_{m=1}^M \int_{\Re} h^m(x) q_m^\ast(x,t)\ud x$;
\item $\widehat{h_t^m} :=\E[h(X_t,\theta_t)|\theta_t = e_m,\UZ] = \int_{\Re} h^m(x) \rho^\ast_m(x,t) \ud x$.
\end{romannum}
Note that $\hat{h}_t = \sum_{m=1}^M \mu_t^m \widehat{h_t^m}$.

\medskip

The following theorem describes the evolution of above-mentioned density functions. A short proof is included in Appendix~\ref{apdx:jointpm}. 

\begin{theorem}[See also Theorem 1 in~\cite{Blom_cdc12}]
\label{thm:jointpm_evol}
Consider the hybrid system~\eqref{eqn:process_model}\,-\,\eqref{eqn:observ_model}:
\begin{romannum}
\item The joint conditional distribution of $(X_t, \theta_t)^T$ satisfies:
\begin{align}
\ud q^\ast = \mathcal{L}^\dagger (q^\ast) \ud t + Q^T q^\ast \ud t + (H_t - \hat{h}_t I)(\ud Z_t - \hat{h}_t\ud t)q^\ast,\label{eqn:kushner_q}
\end{align}
where $\Ldag = \text{diag}\{\Ldag_m\}$, $H_t = \text{diag}\{h^m\}$, $I$ is an $M\times M$ identity matrix and
\begin{equation}
\Ldag_m q_m^\ast := -\frac{\partial}{\partial x}(a^m q_m^\ast) + \frac{1}{2}(\sigma^m)^2 \frac{\partial^2}{\partial x^2} q_m^\ast.\nonumber
\end{equation}
\item The conditional distribution of $\theta_t$ satisfies:
\begin{equation}
\ud \mu_t = Q^T \mu_t \ud t +  (\widehat{H}_t - \hat{h}_t I) (\ud Z_t - \hat{h}_t \ud t)\mu_t ,\label{eqn:kushner_mu}
\end{equation}
where
$\widehat{H}_t = \text{diag}\{\widehat{h_t^m}\}$.
\item The conditional distribution of $X_t$ satisfies:
\begin{equation}
\ud p^\ast = \sum_{m=1}^M \Ldag_m(q_m^\ast) \ud t + \sum_{m=1}^M (h^m - \hat{h}_t) (\ud Z_t - \hat{h}_t \ud t) q_m^\ast.\label{eqn:kushner_p}
\end{equation}
\item The conditional distribution of $X_t$ given $\theta_t = e_m$ satisfies:
\begin{align}
\ud \rho_m^\ast = &\Ldag_m \rho_m^\ast  \ud t  + \frac{1}{\mu_t^m}\sum_{l=1}^M q_{lm}\mu_t^l (\rho_l^\ast  - \rho_m^\ast ) \ud t  \nonumber\\
&+ (h^m - \widehat{h_t^m}) (\ud Z_t - \widehat{h_t^m}\ud t)\rho_m^\ast ,\quad m = 1,\hdots, M\label{eqn:kushner_rho}
\end{align}

\end{romannum}
\qed
\end{theorem}

\section{IMM-Feedback Particle Filter}
\label{sec:Ifilter}

The IMM-FPF is comprised of $M$ parallel feedback particle filters: The model for the $\mth$ particle filter is given by,
\begin{equation}
\ud X_t^{i;m} = a^m(X_t^{i;m}) \ud t + \sigma^m \ud B_t^{i;m} + \ud U_t^{i;m},\label{eqn:ifilter_control}
\end{equation}
where $X_t^{i;m}\in \Re$ is the state for the $\ith$ particle at time $t$, $U_t^{i;m}$ is its control input, $\{B_t^{i;m}\}_{i=1}^N$ are mutually independent standard Wiener processes and $N$ denotes the number of particles. We assume that the initial conditions $\{X^{i;m}_0\}_{i=1}^N$  are
i.i.d., independent of $\{B^{i;m}_t\}$, and drawn from the initial
distribution $\rho_m^*(x,0)$ of $X_0$.  Both  $\{B^{i;m}_t\}$ and
$\{X^{i;m}_0\}$ are assumed to be independent of $X_t,Z_t$.
Certain additional assumptions are made regarding admissible
forms of control input (see~\cite{taoyang_cdc12}).

\begin{remark}
The motivation for choosing the parallel structure comes from the conventional IMM filter, which is comprised of $M$ parallel Kalman filters, one for each maneuvering mode $m\in \{1,\hdots,M\}$.
\qed
\end{remark}

There are two types of conditional distributions of
interest in our analysis:
\begin{romannum}
\item $\rho_m^*(x,t)$: Defines the conditional dist.\ of
    $X_t$ given $\theta_t = e_m$ and $\UZ$, see~\eqref{eqn:def_rho_m}.
\item $\rho_m(x,t)$:  Defines the conditional dist. of $X_t^{i;m}$ given $\UZ$:
    \begin{equation}
    \int_{x\in A} \rho_m(x,t) \ud x = \P\{X_t^{i;m}\in A|\UZ\},\quad \forall A\in \mathcal{B}(\Re).\nonumber
    \end{equation}

\end{romannum}

The control problem is to choose the control inputs $\{U^{i;m}_t\}_{m=1}^M$ so
that, $\rho_m$ approximates $\rho_m^*$ for each $m = 1, \hdots, M$. Consequently the empirical
distribution of the particles approximates $\rho_m^*$ for large
number of particles~\cite{taoyang_acc12}.

The main result of this section is to describe an explicit
formula for the optimal control input, and demonstrate that
under general conditions we obtain an exact match:  $\rho_m=\rho_m^*$, under optimal control. The optimally controlled dynamics of the
$i^{\text{th}}$ particle in the $\mth$ FPF have the following Stratonovich form,
\begin{align}
\ud X_t^{i;m} = a^m(X_t^{i;m})\ud t &+ \sigma^m \ud B_t^{i;m} + \v^m(X_t^{i;m},t) \circ \ud I_t^{i;m} \nonumber\\
&+ u^m(X_t^{i;m}, X_t^{i;-m}) \ud t,\label{eqn:IMMFPF_ind}
\end{align}
where $X_t^{i;-m} = \{X_t^{i;l}\}_{l \neq m}$ and $\Inov^{i;m}_t$ is the modified form of
$\emph{innovation process}$,
\begin{equation}
\ud I_t^{i;m} := \ud Z_t - \frac{1}{2}\left[h^m(X_t^{i;m}) + \widehat{h_t^m}\right]\ud t,\label{eqn:IMMFPF_inov}
\end{equation}
where $\widehat{h_t^m}:= \int_{\Re} h^m(x) \rho_m(x,t)\ud x$. The gain function $\v^m$ is obtained as a solution of an Euler-Lagrange boundary value problem (E-L BVP):
\begin{equation}
\frac{\partial (\rho_m \v^m)}{\partial x} = -(h^m-\widehat{h_t^m})\rho_m,\label{eqn:pasthatm_bvp}
\end{equation}
with boundary condition $\lim_{x\rightarrow \pm \infty} \rho_m(x,t) \v^m(x,t) = 0$.

The interaction between filters arises as a result of the control term $u^m$. It is obtained by solving the following BVP:
\begin{equation}
\frac{\partial (\rho_m u^m)}{\partial x} = \sum_{l=1}^M c_{lm} (\rho_m - \rho_l),\label{eqn:bvp_utm}
\end{equation}
again with boundary condition $\lim_{x\rightarrow \pm \infty} \rho_m(x,t) u^m(x,t) = 0$ and $c_{lm} := \frac{q_{lm}\mu_t^l}{\mu_t^m}$.

Recall that the evolution of  $\rho_m^*(x,t)$ is described by the modified Kushner-Stratonovich (K-S) equation~\eqref{eqn:kushner_rho}. The evolution of $\rho_m$ is given by a forward
Kolmogorov operator (derived in Appendix~\ref{apdx:pf_consistency}).

The following theorem shows that the evolution equations
for $\rho_m$ and $\rho_m^\ast$ are identical. The proof appears in
Appendix~\ref{apdx:pf_consistency}.

\begin{theorem}
\label{thm:rhom_consistency}
Consider the two distributions $\rho_m$ and $\rho_m^\ast$. 
Suppose that, for $m = 1,\hdots, M$, the gain function $\v^m(x,t)$ and the control term $u^m$ is obtained according to~\eqref{eqn:pasthatm_bvp} and~\eqref{eqn:bvp_utm}, respectively. Then provided $\rho_m(x,0) = \rho_m^\ast(x,0)$, we have for all $t > 0$, $\rho_m(x,t) = \rho_m^\ast(x,t)$.\qed
\end{theorem}

\medskip

In a numerical implementation, one also needs to estimate $\mu_t^m$, which is done by using the same finite-dimensional filter, as in~\eqref{eqn:kushner_mu}:
\begin{equation}
\qquad \ud \mu_t^m = \sum_{l=1}^M q_{lm}\mu_t^l \ud t + (\widehat{h_t^m} - \hat{h}_t) (\ud Z_t - \hat{h}_t \ud t)\mu_t^m,\label{eqn:mu_m_evol}
\end{equation}
where $\hat{h}_t = \sum_{m=1}^M \mu_t^m \widehat{h_t^m}$ and $\widehat{h_t^m} \approx \frac{1}{N} \sum_{i=1}^N h^m(X_t^{i;m})$ are approximated with particles.
\begin{remark}
The mode association probability filter~\eqref{eqn:mu_m_evol} can also be derived by considering a continuous-time limit starting from the continuous-discrete time filter that appears in the classic IMM filter literature~\cite{BarshalomLiKirubarajan01}.
This proof appears
in Appendix~\ref{apdx:discrete_assoc}. The alternate proof
is included because it shows that the filter~\eqref{eqn:mu_m_evol} is in fact the continuous-time nonlinear counterpart of the algorithm that is used to obtain association probability in the classical IMM literature. The proof also suggests alternate discrete-time algorithms for evaluating
  association probabilities in simulations and experiments, where
  observations are made at discrete sampling times. 
  \qed
\end{remark}

\smallskip

Define $p(x,t) := \sum_{m=1}^M \mu_t^m \rho_m(x,t)$ where $\mu_t^m$ is defined in~\eqref{eqn:def_mu_m}. The following corollary shows that $p(x,t)$ and $p^\ast(x,t)$ are identical. Its proof is straightforward, by using the definitions, and is thus omitted.
\begin{corollary}
Consider the two distribution $p(x,t)$ and $p^\ast(x,t)$. Suppose conditions in~\Thm{thm:rhom_consistency} apply, and $\mu_t^m$ is obtained using~\eqref{eqn:mu_m_evol}, then provided $p(x,0) = p^\ast(x,0)$, we have for all $t > 0$, $p(x,t) = p^\ast(x,t)$.\qed
\end{corollary}

\subsection{Algorithm}
The main difficulty is to obtain a solution of the BVP~\eqref{eqn:pasthatm_bvp} and~\eqref{eqn:bvp_utm} at each time step. A Galerkin algorithm for the same appears in our earlier papers~\cite{taoyang_cdc12},\cite{adamtilton_acc13}. One particular approximation of the solution, referred to as the {\em constant gain approximation} is given by:
\begin{align}
\v^m &\approx \frac{1}{N}\sum_{i=1}^N (h^m(X_t^{i;m}) - \widehat{h_t^m})X_t^{i;m} ,\label{eqn:bvpm_approx}\\
u^m &\approx \sum_{l=1}^M c_{lm} \frac{1}{N}\left(\sum_{i=1}^N X_t^{i;l} - \sum_{i=1}^N X_t^{i;m}\right).\label{eqn:utm_const}
\end{align}
The derivation of the constant approximation~\eqref{eqn:bvpm_approx}-\eqref{eqn:utm_const} appears in Appendix~\ref{apdx:deri_const_approx}.

Apart from the gain function, the algorithm requires approximation of $\hat{h}_t$ and $\widehat{h_t^m}$. These are obtained in terms of particles as:
\begin{align}
\widehat{h_t^m} \approx \frac{1}{N} \sum_{i=1}^N h^m(X_t^{i;m}),\qquad \hat{h}_t = \sum_{m=1}^M \mu_t^m \widehat{h_t^m}.\nonumber
\end{align}
For simulating the IMM-FPF, we use an Euler-discretization method. The resulting discrete-time algorithm appears in~Algo.1. At each time step, the algorithm requires computation of the gain function, that is obtained using~\eqref{eqn:bvpm_approx}-\eqref{eqn:utm_const}.

\begin{algorithm}
\label{algo:ifilter}
\caption{IMM-FPF for SHS}
\begin{algorithmic}[1]
\STATE {INITIALIZATION}
\FOR{$m=1$ to $M$}
\STATE $\mu_0^m = \frac{1}{M}$.
\FOR{$i=1$ to $N$}
\STATE Sample $X_0^{i;m}$ from $p^*(\cdot,0)$.
\ENDFOR
\ENDFOR
\STATE $p^N(x,0) = \frac{1}{N}\sum_{i=1}^N \sum_{m=1}^M \mu_0^m \delta_{X_0^{i;m}}(x)$
\end{algorithmic}
\begin{algorithmic}[1]
\STATE {ITERATION} [$t$ to $t + \Delta t$]
\FOR{$m=1$ to $M$}
\STATE Calculate $\widehat{h_t^m} \approx \frac{1}{N} \sum_{i=1}^N h^m(X_t^{i;m})$.
\ENDFOR
\STATE Calculate $\hat{h}_t = \sum_{m=1}^M \mu_t^m \widehat{h_t^m}$.
\FOR{$m=1$ to $M$}
\FOR{$i =1$ to $N$}
\STATE Generate a sample, $\Delta V$, from $N(0,1)$
\STATE Calculate $\Delta I_t^{i;m} = \Delta Z_t - \frac{1}{2}\left[h^m(X_t^{i;m}) + \widehat{h_t^m}\right]\Delta t$
\STATE Calculate the gain function $\v^m$ (e.g., by using~\eqref{eqn:bvpm_approx})
\STATE Calculate the control term $u^m$ (e.g., by using~\eqref{eqn:utm_const}).
\STATE $\Delta X_t^{i;m} = a^m(X_t^{i;m}) \Delta t + \sigma^m \sqrt{\Delta t} \Delta V + \v^m \Delta I_t^{i;m} + u^m \Delta t$
\STATE $X_{t + \Delta t}^{i;m} = X_t^{i;m} + \Delta X_t^{i;m}$
\ENDFOR
\STATE $\mu_{t+\Delta t}^m = \mu_t^m + \sum_{l=1}^M q_{lm} \mu_t^l \Delta t + (\widehat{h_t^m} - \hat{h}_t) (\Delta Z_t - \hat{h}_t\Delta t) \mu_t^m$.
\ENDFOR
\STATE $t = t + \Delta t$
\STATE $p^N(x,t) = \frac{1}{N}\sum_{i=1}^N \sum_{m=1}^M \mu_t^m \delta_{X_t^{i;m}}(x)$
\STATE $\hat{X}_t := \E[X_t|\UZ] \approx \frac{1}{N}\sum_{m=1}^M \sum_{i=1}^N \mu_t^m X_t^{i;m}$.
\end{algorithmic}
\end{algorithm}

\begin{figure}
\includegraphics[scale = 0.5]{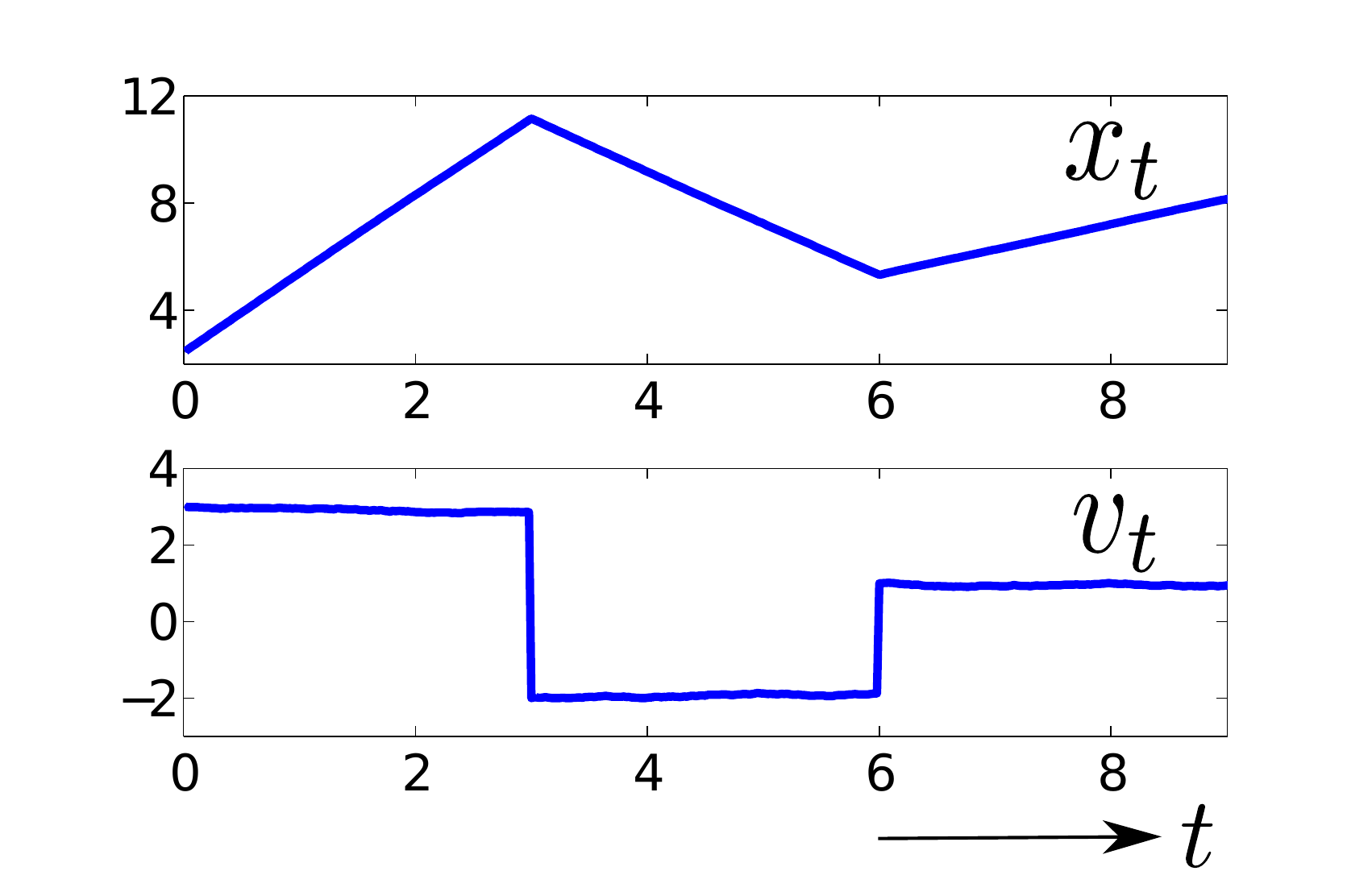}
\caption{Simulation results for a single trial (from top to bottom): (a) Sample
   position path $x_t$; (b) Sample velocity path $v_t$.}
\label{fig:true_traj}
\end{figure}

\section{Numerics}
\label{sec:numerics}

\begin{figure*}
\centering
\includegraphics[scale = 0.5]{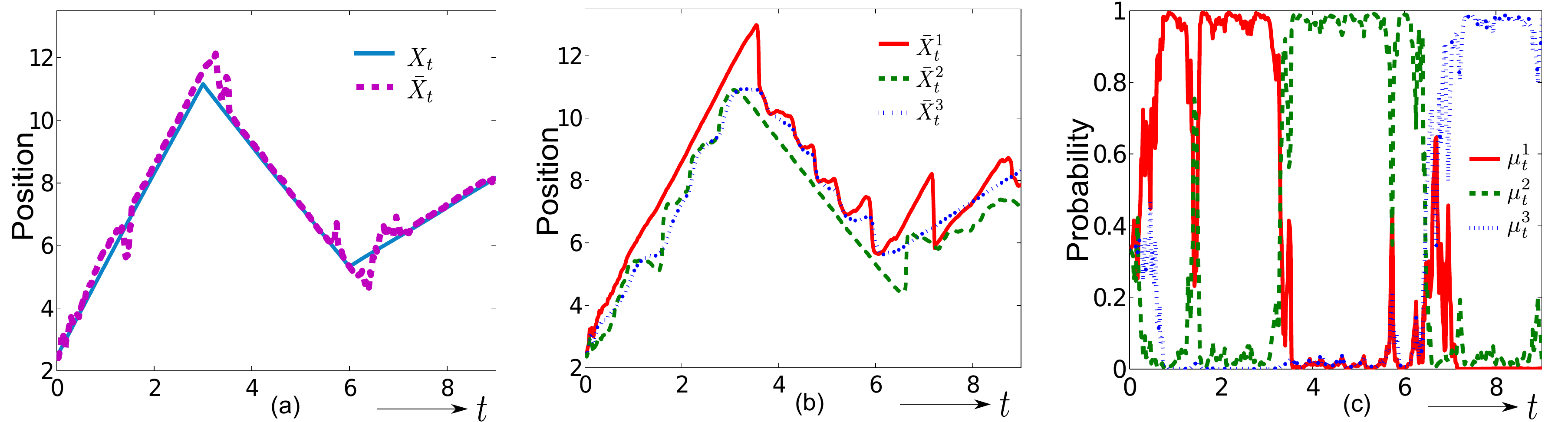}
\caption{Maneuvering target tracking using IMM-FPF: (a) Comparison of IMM-FPF estimated mean $\bar{X}_t$ with target trajectory $X_t$. (b) Plot of mean trajectories of individual modes. (c) Plot of mode association probability.}
\label{fig:ifilter}
\end{figure*}

\subsection{Maneuvering target tracking with bearing measurements}
We consider a target tracking problem where the target dynamics evolve according to a white noise acceleration model:
\begin{align}
\ud X_t &= \begin{bmatrix}
                    0 &  1 \\
                    0 &  0
  \end{bmatrix} X_t \ud t + \sigma_B \ud B_t,\label{eqn:num_dyn}\\
\ud Z_t &= h(X_t) \ud t + \sigma_W \ud W_t,\label{eqn:num_obs}
\end{align}
where $X_t = (x_t,v_t)$ denotes the state vector comprising of position and
velocity coordinates at time $t$, $Z_t$ is the observation
process, and $\{B_t\},\{W_t\}$ are mutually independent standard
Wiener processes. 
We consider a bearing-only measurement with:
\begin{equation}
h(x,v) = \arctan \left(\frac{x}{L}\right),\label{eqn:bear_meas}
\end{equation}
where $L$ is a constant.

In the target trajectory simulation, the following parameter values are used: $\sigma_B = [0, 0.05]$, $\sigma_W =  0.015$, $L = 10$ and initial condition $X_0 =(x_0,v_0) = (2.5, 3)$. The total simulation time is $T = 9$ and the time step for numerical discretization is $\Delta t = 0.02$. At $T_1 = 3$ and $T_2 = 6$, the target instantaneously changes its velocity to $v = -2$ and $v = 1$ respectively. The resulting trajectory is depicted in Figure~\ref{fig:true_traj}. At each discrete time step, a bearing measurement is obtained according to~\eqref{eqn:bear_meas}. The target is maneuvering in the sense that its velocity switches between three different values $\{3,-2,1\}$.


\subsection{Tracking using IMM-FPF}

We assume an interacting multiple model architecture as follows:
\begin{romannum}
\item There are three possible target dynamic modes:
\begin{align}
\theta_t = 1:\qquad \ud X_t &= 3 \; \ud t + \sigma_B \ud B_t,\nonumber\\
\theta_t = 2:\qquad \ud X_t &= -2 \; \ud t + \sigma_B \ud B_t,\nonumber\\
\theta_t = 3:\qquad \ud X_t &= 1 \; \ud t + \sigma_B \ud B_t.\nonumber
\end{align}
\item $\theta_t$ evolves as a continuous-time Markov chain with transition rate matrix $Q$.
\item Observation process is modeled the same as in~\eqref{eqn:num_obs}-\eqref{eqn:bear_meas}.
\end{romannum}

In the simulation results described next, we use the following parameter values: $\sigma_B = 0.05$, $\sigma_W = 0.015$ and initial condition $X_0 = 2.5$. The total simulation time is $T = 9$ and time step $\Delta t = 0.02$. The transition rate matrix is,
\begin{equation*}
Q = \begin{bmatrix}
                    -0.1 &  0.1 & 0 \\
                    0.05 &  -0.1& 0.05\\
                    0 & 0.1& -0.1
  \end{bmatrix}
\end{equation*}
The prior mode association probability $\mu_t = (\mu_t^1,\mu_t^2,\mu_t^3)$ at time $t = 0$ is assumed to be $\mu_0 = (1/3,1/3,1/3)$.

Figure~\ref{fig:ifilter}(a) depicts the result of a single simulation: The estimated mean trajectories, depicted as dashed lines, are obtained using the IMM-FPF algorithm described in~\Sec{sec:Ifilter}. Figure~\ref{fig:ifilter}(b) depicts the mean trajectories of individual modes. Figure~\ref{fig:ifilter}(c) depicts the evolution of association probability $(\mu_t^1,\mu_t^2,\mu_t^3)$ during the simulation run. We see that trajectory probability converges to the correct mode during the three maneuvering periods. For the filter simulation, $N = 500$ particles are used for each mode. 



\section{Conclusion and Future Work}
\label{sec:conclustion}

In this paper, we introduced a feedback particle filter-based algorithm for the continuous-time SHS estimation problem. The proposed algorithm is shown to be the nonlinear generalization of the conventional Kalman filter-based IMM algorithm. A numerical example for a maneuvering target tracking problem is presented to illustrate the use of IMM-FPF.

The ongoing research concerns the following two topics:
\begin{romannum}
\item Comparison of the IMM-FPF against the basic particle filter (PF) and IMM-PF, with respect to estimation performance and computational burden.
\item Investigation of alternate FPF-based algorithms for SHS estimation. Of particular interest is the filter architecture where particles evolve on the joint state space (analogous to case (i) in~\Sec{sec:intro}).
\end{romannum}

\appendix

\subsection{Derivation of the exact filtering equations~\eqref{eqn:kushner_q}-\eqref{eqn:kushner_rho}}
\label{apdx:jointpm}

\head{\bf {Derivation of~\eqref{eqn:kushner_q}:}}
For each fixed mode $\theta_t = e_m \in \{e_1,e_2,\hdots,e_M\}$, the state process $X_t$ is a Markov process with Kolmogorov forward operator $\mathcal{L}^{\dagger}_m$. Therefore, the joint process $(X_t,\theta_t)^T$ is a Markov process with generator $\mathcal{L}^\dagger + Q$ where $\mathcal{L}^\dagger := \text{diag} \{\mathcal{L}^\dagger_m\}$. Defining $\int_A \pi^\ast_m(x,t) \ud x := \P\{X_t \in A, \theta_t = e_m\}$ and $\pi^\ast(x,t) := (\pi_1^\ast(x,t),\cdots,\pi^\ast_M(x,t))^T$ we have that
\begin{equation}
\frac{\partial \pi^\ast}{\partial t} (x,t) = \mathcal{L}^\dagger \pi^\ast(x,t) + Q^T \pi^\ast(x,t).\nonumber
\end{equation}

Recall that the posterior distribution was defined by $\int_A q^\ast_m(x,t) \ud x := \P\{[ X_t \in A, \theta_t = e_m]|\UZ\}$. By applying the fundamental filtering theorem for Gaussian observations (see~\cite{Lipster77}) to~\eqref{eqn:observ_model} we have:
\begin{align}
\ud q^\ast_m = &\Ldag_m(q^\ast_m) \ud t + \sum_{l=1}^M q_{lm} q^\ast_l \ud t + (h^m-\hat{h}_t) (\ud Z_t - \hat{h}_t \ud t) q^\ast_m,\nonumber
\end{align}
where $
\hat{h}_t := \sum_{m=1}^M \int_{\Re}  h^m(x) q^\ast_m(x,t) \ud x
$.\qed

\smallskip

\head{\bf {Derivation of~\eqref{eqn:kushner_mu} and~\eqref{eqn:kushner_p}:}}
By definition, we have,
\begin{align}
p^\ast (x,t) &= \sum_{m=1}^M q_m^\ast(x,t),\label{eqn:cond_em}\\
\mu_t^m &= \int_\Re q_m^\ast(x,t) \ud x.\label{eqn:cond_mu}
\end{align}
Taking derivatives on both sides of~\eqref{eqn:cond_em} and~\eqref{eqn:cond_mu} gives the desired result.\qed

\smallskip

\head{\bf {Derivation of~\eqref{eqn:kushner_rho}:}}
By definition $q^\ast_m = \rho_m^\ast \mu_t^m$. Applying It$\hat{\text{o}}$'s differentiation rule we have:
\begin{equation}
\ud \rho_m^\ast = \frac{\ud q^\ast_m}{\mu_t^m} + q^\ast_m \ud \left(\frac{1}{\mu_t^m}\right) + \ud q^\ast_m \ud \left(\frac{1}{\mu_t^m}\right),\label{eqn:chain_rule}
\end{equation}
where $\ud \left(\frac{1}{\mu_t^m}\right) = -\frac{\ud \mu_t^m}{(\mu_t^m)^2} + \frac{(\ud \mu_t^m)^2}{(\mu_t^m)^3}$. Substituting~\eqref{eqn:kushner_q} and~\eqref{eqn:kushner_mu} into~\eqref{eqn:chain_rule} we obtain the desired result.\qed

\subsection{Proof of consistency for IMM-FPF}
\label{apdx:pf_consistency}
We express the feedback particle filter~\eqref{eqn:IMMFPF_ind} as:
\begin{align}
\ud X_t^{i;m} = a^m(X_t^{i;m})\ud t &+ \sigma^m \ud B_t^{i;m} + \v^m(X_t^{i;m},t)\ud Z_t + \tilde{u}(X_t^{i;m},t) \ud t \nonumber\\
 &+ u^m(X_t^{i;m}, X_t^{i;-m}) \ud t,\nonumber
\end{align}
where
\begin{equation}
\tilde{u}(x,t) = -\frac{1}{2}\v^m(x,t)(h^m(x) + \widehat{h_t^m}) + \Omega(x,t),\label{eqn:def_u_tilde}
\end{equation}
and $\Omega := \frac{1}{2}\v^m (\v^m)'$ is the Wong-Zakai correction term for~\eqref{eqn:IMMFPF_ind}.
The evolution equation for $\rho_m$ is given by:
\begin{equation}
\begin{aligned}
\ud \rho_m = &\clL_m^\dagger \rho_m\ud t  - \frac{\partial}{\partial x} (\rho_m\v^m) \ud Z_t - \frac{\partial}{\partial x}(\rho_m \tilde{u}) \ud t
\\
&  - \frac{\partial}{\partial x}(\rho_m u^m)\ud t +
\frac{1}{2}\frac{\partial^2}{\partial x^2} \left( \rho_m (\v^m)^2 \right) \ud t.
  \end{aligned}
\label{eqn:mod_FPK}
\end{equation}
The derivation of this equation is similar to the basic FPF case (see Proposition 2 in~\cite{taoyang_cdc11}) and thus omitted here.

It is only necessary to show that with the choice of $\{\v^m,u^m\}$ according to~\eqref{eqn:pasthatm_bvp}-\eqref{eqn:bvp_utm},  we have $\ud \rho_m(x,t) = \ud \rho_m^*(x,t)$, for all $x$ and
$t$,  in the sense that they are defined by identical
stochastic differential equations.   Recall $\ud \rho_m^*$ is
defined according to the modified K-S equation~\eqref{eqn:kushner_rho},
and $\ud \rho_m$ according to the forward
equation~\eqref{eqn:mod_FPK}.

If $\v^m$ solves the E-L BVP~\eqref{eqn:pasthatm_bvp} then we have:
\begin{equation}
\frac{\partial}{\partial x} (\rho_m \v^m) = -(h^m - \widehat{h_t^m}) \rho_m. \label{eqn:recall_utm}
\end{equation}
On multiplying both sides of~\eqref{eqn:def_u_tilde} by $-\rho_m$, we have:
\begin{align*}
- \rho_m  \tilde{u}&= \frac{1}{2}(h^m - \widehat{h_t^m}) \rho_m \v^m - \frac{1}{2}(\rho_m \v^m) \frac{\partial \v^m}{\partial x} + \widehat{h_t^m}\rho_m \v^m\\
					 &= -\frac{1}{2}\frac{\partial (\rho_m \v^m)}{\partial x} \v^m - \frac{1}{2}(\rho_m \v^m) \frac{\partial \v^m}{\partial x} + \widehat{h_t^m} \rho_m \v^m\\
					 &= -\frac{1}{2}\frac{\partial }{\partial x}\left(\rho_m (\v^m)^2\right) + \widehat{h_t^m} \rho_m \v^m
\end{align*}
where we used~\eqref{eqn:recall_utm} to obtain the second equality. Differentiating once with respect to $x$ and using~\eqref{eqn:recall_utm} once again,
\begin{equation}
- \frac{\partial}{\partial x} (\rho_m  \tilde{u}) + \frac{1}{2} \frac{\partial^2}{\partial x^2}\left(\rho_m (\v^m)^2\right) = - \widehat{h_t^m}(h^m - \widehat{h_t^m})\rho_m. \label{eqn:deriv_1}
\end{equation}

Substituting~\eqref{eqn:pasthatm_bvp}, \eqref{eqn:bvp_utm} and~\eqref{eqn:deriv_1} to~\eqref{eqn:mod_FPK} and after some simplifications, we obtain:
\begin{align}
\ud \phatm = \Ldag_m \phatm \ud t  &+ (h^m - \widehat{h_t^m}) (\ud Z_t - \widehat{h_t^m}\ud t)\phatm \nonumber\\
&+ \frac{1}{\mu_t^m}\sum_{l=1}^M q_{lm}\mu_t^l (\rho_l - \rho_m) \ud t.\nonumber
\end{align}
This is precisely the SDE~\eqref{eqn:kushner_rho}, as desired.
\qed

\subsection{Alternate Derivation of~\eqref{eqn:mu_m_evol} }
\label{apdx:discrete_assoc}

The aim of this section is to derive, formally, the update part
of the continuous time
filter~\eqref{eqn:mu_m_evol} by taking a
continuous time limit of the discrete-time algorithm for
evaluation of association probability. The procedure for taking
the limit is similar to Sec~$6.8$ in~\cite{jazwinski70} for
derivation of the K-S equation.

At time $t$, we have $M$ possible modes for the SHS. The discrete-time filter for mode association probability is obtained by using Bayes' rule (see~\cite{BarshalomLiKirubarajan01}):
\begin{equation}
\Prob\{\theta_t = e_m|\UZ,\Delta Z_t\} = \frac{\Prob\{\Delta Z_t|\theta_t = e_m\}\Prob\{\theta_t = e_m|\UZ\}}{\sum_{l=1}^M\Prob\{\Delta Z_t|\theta_t = e_l\}\Prob\{\theta_t = e_l|\UZ\}}.\label{eqn:discret_bayes}
\end{equation}

Rewrite:
\begin{align}
\P\{\Delta Z_t|\theta_t =e_m\} &= \int \P\{\Delta Z_t|\theta_t=e_m,X_t=x\}\rho_m(x,t)\ud x\nonumber\\
                          &= L_m(\Delta Z_t).\label{eqn:Z_t_m}
\end{align}
where $L_m(\Delta Z_t) := \frac{1}{\sqrt{2\pi\Delta
t}}\int_\mathbb{R} \exp\left[-\frac{(\Delta Z_t - h^m(x)\Delta
t)^2}{2  \Delta t}\right]\rho_m(x,t)\ud x$.

Now, recall $\mu_t^m = \P\{\theta_t=e_m|\UZ\}$, the increment in
the measurement update step (see Sec~$6.8$
in~\cite{jazwinski70}) is given by
\begin{equation}
\Delta \mu_t^m:= \P\{\theta_t=e_m|\UZ,\Delta Z_t\} - \P\{\theta_t=e_m|\UZ\}.\label{eqn:def_dbeta}
\end{equation}
Using~\eqref{eqn:discret_bayes} and~\eqref{eqn:def_dbeta}, we
have:
\begin{equation}
\Delta \mu_t^m = E^m(\Delta t,\Delta Z_t)\mu_t^m  - \mu_t^m ,\label{eqn:cal_dbeta}
\end{equation}
where
\begin{equation}
E^m(\Delta t,\Delta Z_t) = \frac{\P\{\theta_t = e_m|\UZ,\Delta Z_t\}}{\P\{\theta_t = e_m|\UZ\}}.\label{eqn:E_def}
\end{equation}
We expand $E^m(\Delta t, \Delta Z_t)$ as a multivariate
series about $(0,0)$:
\begin{align}
E^m(\Delta t, \Delta Z_t) &= E^m(0,0) + E^m_{\Delta t}(0,0)\Delta t + E^m_{\Delta Z_t} (0,0)\Delta Z_t\nonumber\\
                                &+ \frac{1}{2} E^m_{\Delta Z_t,\Delta Z_t}(0,0)\ud Z^2_t + o(\Delta t).\label{eqn:multi_series}
\end{align}

By direct evaluation, we obtain:
\begin{align}
&E^m(0,0) = 1,\nonumber\\
&E^m_{\Delta t}(0,0) = \frac{1}{2}\left(\widehat{h_t^2} - \widehat{(h_t^m)^2}\right),\nonumber\\
&E^m_{\Delta Z_t} (0,0) = \widehat{h_t^m} - \hat{h}_t\nonumber\\
&E^m_{\Delta Z_t,\Delta Z_t}(0,0) = \widehat{(h_t^m)^2}-2\widehat{h_t^m}\hat{h}_t + 2 \hat{h}_t^2 - \widehat{h_t^2} \nonumber
\end{align}
where $\widehat{(h_t^m)^2} := \int_{\Re} (h^m(x))^2 \rho_m(x,t) \ud x$ and
$\widehat{h_t^2}:=\sum_{m=1}^M \mu_t^m \widehat{(h_t^m)^2}$.

\medskip

By using It$\hat{\text{o}}$'s rules,
\begin{equation}
\E[\Delta Z_t \Delta Z_t] = \Delta t.\nonumber
\end{equation}
This gives
\begin{align}
E^m(\Delta t,\Delta Z_t) &=  1 + (\widehat{h_t^m} - \hat{h}_t) (\Delta Z_t - \hat{h}_t\Delta t)\label{eqn:Em_approx}
\end{align}
Substituting~\eqref{eqn:Em_approx} to~\eqref{eqn:cal_dbeta} we obtain the expression for $\Delta \mu_t^m$, which equals the measurement update part of the continuous-time
filter.
\begin{remark}
During a discrete-time implementation, one can
use~\eqref{eqn:discret_bayes}-\eqref{eqn:Z_t_m} to obtain
association probability. In~\eqref{eqn:discret_bayes}, $L_m(\ud
Z_t)$ is approximated by using particles:
\begin{equation}
L_m(\Delta Z_t) \approx \frac{1}{N}\frac{1}{\sqrt{2\pi\Delta t}} \sum_{i=1}^N \exp\left[-\frac{(\Delta Z_t - h^m(X_t^{i;m})\Delta t)^2}{2\Delta t}\right].\nonumber
\end{equation}
\end{remark}

\begin{figure}
\centering
\includegraphics[scale = 0.5]{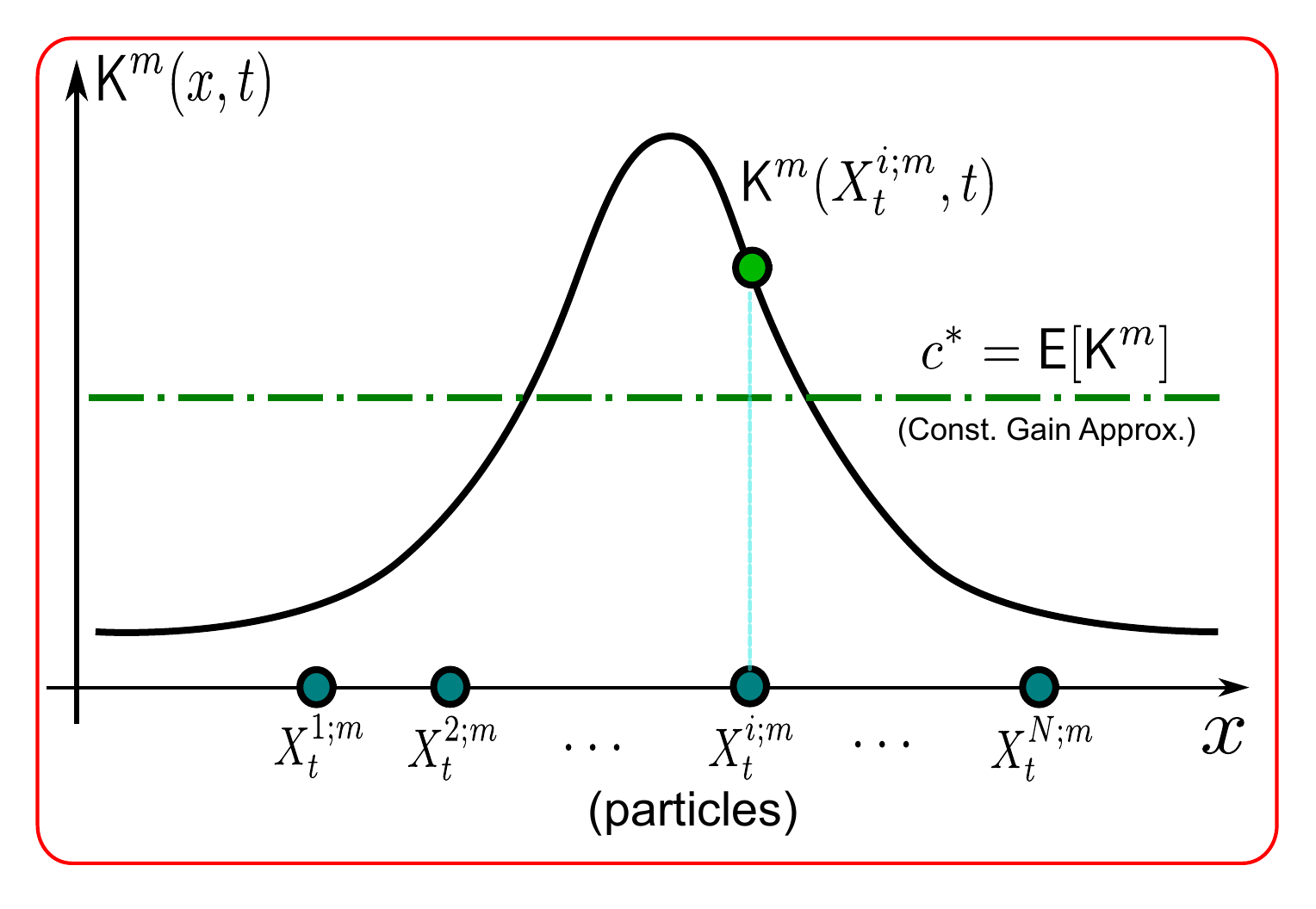}
\caption{Approximating nonlinear $\v^m$ by its expected value $\E[\v^m]$.}
\label{fig:gain_func}
\end{figure}

\subsection{Derivation of constant approximation~\eqref{eqn:bvpm_approx}-\eqref{eqn:utm_const}}
\label{apdx:deri_const_approx}

In this section, we provide a justification for~\eqref{eqn:bvpm_approx}-\eqref{eqn:utm_const}. Recall that at each fixed time step $t$, $\v^m(x,t)$ is obtained by solving the BVP~\eqref{eqn:pasthatm_bvp}:
\begin{equation*}
\frac{\partial (\rho_m \v^m)}{\partial x} = -(h^m-\widehat{h_t^m})\rho_m.
\end{equation*}
A function $\v^m$ is said to be a weak solution of the BVP~\eqref{eqn:pasthatm_bvp} if
\begin{equation}
\E\left[\v^m \frac{\partial \psi}{\partial x}\right] = \E[(h^m-\widehat{h_t^m})\psi]\label{eqn:weak_form}
\end{equation}
holds for all $\psi \in H^1(\Re;\rho_m)$ where $\E[\cdot] := \int_\Re \cdot \rho_m(x,t)\ud x$ and $H^1$ is a certain Sobolev space (see~\cite{taoyang_cdc12}). The existence-uniqueness results for the weak solution of~\eqref{eqn:weak_form} also appear in~\cite{taoyang_cdc12}.

In general, the weak solution $\v^m(\cdot, t)$ of the BVP~\eqref{eqn:weak_form} is some nonlinear scalar-valued function of the state (see~\Fig{fig:gain_func}). The idea behind the {\em {constant gain approximation}} is to find a single constant $c^\ast \in \Re$ to approximate this function (see~\Fig{fig:gain_func}). Precisely,
\begin{equation*}
c^\ast = \arg \min_{c \in \Re} \E[(\v^m - c)^2].
\end{equation*}
By using a standard sum of square argument, we have
\[
c^\ast = \E[\v^m].
\]
Even though $\v^m$ is unknown, the constant $c^\ast$ can be obtained using~\eqref{eqn:weak_form}. Specifically, by substituting $\psi(x) = x$ in~\eqref{eqn:weak_form}:
\begin{align}
\E[\v^m] = \E[(h^m-\widehat{h_t^m})\psi] = \int_\Re (h^m(x) - \widehat{h_t^m})\; x\; \rho_m(x,t)\ud x.\nonumber
\end{align}
In simulations, we approximate the last term using particles:
\begin{equation*}
\E[\v^m] \approx \frac{1}{N}\sum_{i=1}^N \left(h^m(X_t^{i;m}) - \widehat{h_t^m}\right) X_t^{i;m},
\end{equation*} 
which gives~\eqref{eqn:bvpm_approx}. The derivation for~\eqref{eqn:utm_const} follows similarly.\qed


\bibliographystyle{plain}
\bibliography{13CDC_IMMFPF,fpfbib,13ACC_DSFPF,refmtt,ACCJPDAFPF}
\end{document}